\documentclass{ifacconf}

\usepackage{graphicx}      % include this line if your document contains figures
\usepackage{natbib}        % required for bibliography

\usepackage{amsmath,amsfonts}
\usepackage{bigints}
\usepackage{todonotes}

\newcommand{\Swarm}{{\mathcal S}}

\newcommand{\Rconf}{{\mathcal R}}
\newcommand{\bx}{{\bf x}}
\newcommand{\br}{{\bf r}}
\newcommand{\bv}{{\bf v}}
\begin{document}
\begin{frontmatter}

\title{Velocity Field Generation for Density Control of Swarms using Heat Equation and Smoothing Kernels\thanksref{footnoteinfo}} 

\thanks[footnoteinfo]{This research was supported in part by Defense Advanced Research Projects Agency (DARPA) Grant No. D14AP00084, and National Science Foundation Grant No. CNS-1624328.}

\author[First]{Utku Eren} and 
\author[Second]{Beh\c cet\ A\c c\i kme\c se} 
%\author[Third]{Third C. Author}

\address[First]{University of Washington, 
   Seattle, WA 98195 USA \\ (e-mail: ue@uw.edu).}
\address[Second]{University of Washington, 
   Seattle, WA 98195 USA \\ (e-mail: behcet@uw.edu).}
%\address[Third]{Electrical Engineering Department, 
%   Seoul National University, Seoul, Korea, (e-mail: author@snu.ac.kr)}

%===============================================================================
%===============================================================================
\begin{abstract}                % Abstract of not more than 250 words.
	This paper presents a method to control the probability density distribution of a swarm of vehicles via velocity fields. The proposed approach synthesizes smooth velocity fields, which specify a desired velocity as a function of time and position in a  decentralized manner i.e., each agent calculates the desired velocity locally by utilizing the number of agents within a prescribed communication distance. Swarm converges to the desired/commanded density distribution by following the velocity field. The local information  consists, only, of agents' positions  and it is utilized to estimate the density around each agent. Local density estimation  is performed  by  using \emph{kernel density estimation}. Then the local density estimates and the desired density are utilized to generate the velocity field, which propagates the swarm probability density distribution  
via the well-known \emph{heat equation}. 
%More specifically, local desired velocities are calculated such that the  resulting swarm density evolution is determined by the heat equation, which is a partial differential equation that drives the swarm density to the desired density. 
The key advantage of using smooth velocity fields to control swarm density with respect to our earlier Markov chain based probabilistic control methods is that the agents move more smoothly in an organized manner and their relative velocities go to zero as they get closer to each other, i.e., they facilitate conflict/collision avoidance.  The desired density distribution, which is commanded to the local controllers, can be computed by using a Markov chain that propagates  the desired density distribution within prescribed requirements. Along with convergence and stability analysis, the effectiveness of the approach is illustrated via numerical examples. 
\end{abstract}


\end{frontmatter}

%===============================================================================
%===============================================================================
\section{Introduction}

The idea of using a large number of autonomous vehicles has emerged as a new paradigm over the years. Mainly inspired by natural phenomena of collectives (e.g., social insects), this new concept called \emph{swarm robotics} has a potential to facilitate many new applications. Employing many low-cost agents rather than more capable yet expensive few vehicles has benefits on redundancy, reconfigurability and parallelism. As much as its benefits, swarm robotics introduces many extreme challenges in guidance and control field. For instance, algorithmic scalability is essential to enable control of swarms, thus decentralization is crucial as each agent comes with a computational resource. Also, robustness to dynamic changes in both swarm (e.g., agent addition/subtraction) and the environment by harnessing only local information while staying alerted for conflict/collision avoidance are basis yet not easily achievable requirements of many swarm applications. 

The control problem of swarms have been investigated through various approaches such as leader-follower \cite{desai2001modeling,yu2010distributed} based solutions, algorithms that takes advantage of graph theory and consensus notion \cite{mesbahi2010graph,chapman2015advection,galbusera2007control}, distributed optimization techniques \cite{zhu2015distributed} and artificial potential field approaches \cite{chaimowicz2005controlling}. Also, decentralized controllers are designed with Smooth Particle Hydrodynamics method \cite{monaghan2005smoothed} in \cite{pimenta2008control} and \cite{pimenta2013swarm} to benefit from the characteristic motion of fluids for manipulating swarm as an incompressible flow to achieve desired swarm behavior. Most of these methods primarily focus on interactions between agents through putting the main emphasis on inter-agent distances or using repulsive forces in formulations, hence, for more high-level tasks (e.g. pattern generation) and complex swarm behaviors, their capacities become limited or over-constrained. 

Density control of swarms started to attract attention due to its capability to capture complex swarm behaviors. The notion of density has both probabilistic and deterministic sense, which enables both surveillance (in the time domain for a low number of agents) and complex pattern generation (in the spatial domain for a large number of agents). Furthermore, when treated suitably, density based control encapsulates most of the agent level interactions, such as collision avoidance. In \cite{zhao2011density} an example of density control for group motion and segregation is presented while still utilizing repulsive forces. The work in \cite{krishnan2016self} and \cite{krishnan2017distributed} presents a density based boundary and configuration control via pseudo-localization algorithm for self-organizing swarms which idles the need for global positioning system. The probabilistic notion of density control is illustrated in \cite{elamvazhuthi2016coverage} for stochastic coverage via diffusing swarm of robots that take local measurements of an underlying scalar field.

Our main contribution in this paper is a deterministic method that synthesizes a velocity field such that, when followed by agents, it drives the swarm of robots to the desired density distribution. The method is decentralized in the sense that, it depends only on the positions of neighboring agents. The novel idea behind the method is the utilization of heat equation which operates on the local density differences from the desired density, thus, resulting controller acts as a partial differential equation based local density feedback controller. Thanks to certain properties of the heat equation, the velocity field diffuses agents in a locally uniform manner to the desired density profile, hence, inter-agent distances are directly imposed by the desired density profile. As a potential expansion of the current work, we also provided an insight on the probabilistic interpretation of the problem as it is mentioned in \cite{hamann2008framework},\cite{elamvazhuthi2016coverage} and \cite{berman2011design}.

The paper is organized as follows. Section II gives the statement of density control problem from both deterministic and probabilistic perspective. Section III gives a background on kernel density estimation for obtaining local density information. Section IV introduces the heat equation based velocity field generation method that utilizes only local density estimations and it provides analysis on stability and convergence. Section V provides numerical examples to illustrate the efficiency of the approach. Finally, conclusions are stated in Section VI.

%===============================================================================
%===============================================================================
\section{Density Control Problem}\label{Sec: Problem Definition}
This section introduces the notion of swarm density distribution and the formulation of swarm density control problem for autonomous agents, from a  probabilistic  perspective. 

\subsection{Swarm Density Distribution}
Consider a swarm of $N$ agents with point mass dynamics that are distributed over the configuration space $\Rconf \subset \mathbb{R}^d$ with continuous boundary $\partial \Rconf$. We define the physical density of the swarm as the fraction of number of agents per unit volume in case of $N \rightarrow \infty$ which can be described at a point $\textbf{x} \in \Rconf$ at time $t \in \mathbb{R}^+$ as follows,

\begin{equation} \label{eqn: trueDens}
	\rho(t,\bx) = \lim\limits_{\epsilon \rightarrow 0}\Bigg[\lim\limits_{N \rightarrow \infty} \frac{1}{N}\bigg(\frac{{n}_{{\mathcal{B}}_\epsilon(\bx)}}{V_{{\mathcal{B}}_\epsilon(\bx)}} \bigg)\Bigg] 
\end{equation}

where ${\mathcal{B}}_\epsilon(\bx)$ is the $\epsilon$-neighborhood of a point $\textbf{x}$ defined as,

\begin{equation*}
	\mathcal{B}_\epsilon(\bx) = \big\{\xi \in \mathbb{R}^d : \hspace{2mm} \| \xi - \bx\| < \epsilon \big\},
\end{equation*}

${n}_{{\mathcal{B}}_\epsilon(\bx)}$ is the total number of agents within ${\mathcal{B}}_\epsilon(\bx)$ and
$V_{{\mathcal{B}}_\epsilon(\bx)}$ is the volume of ${\mathcal{B}}_\epsilon(\bx)$.
Notice that the term $V_{{\mathcal{B}}_\epsilon(\bx)}$ can be moved out from the inner limit in (\ref{eqn: trueDens}) as $V_{{\mathcal{B}}_\epsilon(\bx)}$ is independent of $N$. Using \emph{law of large numbers} \cite{chung2001course} Theorem 5.4.2, the remaining limit expression can be interpreted as the probability of having an agent in the $\epsilon$-neighborhood of a point $\textbf{x} \in \Rconf$ at time $t \in \mathbb{R}^+$,

\begin{equation} \label{eqn: ProbDef}
	\textnormal{prob}\Big( \br(t) \in \mathcal{B}_\epsilon(\bx) \Big) = \mathbb{E}\bigg[\frac{{n}_{{\mathcal{B}}_\epsilon(\bx)}}{N}\bigg] = \lim\limits_{N \rightarrow \infty} \bigg(\frac{{n}_{{\mathcal{B}}_\epsilon(\bx)}}{N} \bigg),
\end{equation}

where $\mathbb{E}[\hspace{0.5mm} \cdot \hspace{0.5mm}]$ is the expected value. Letting $\br(t)$ be the position of an agent,  this probability can also be written in the following form,

\begin{equation} \label{eqn: ProbDens}
	\textnormal{prob}\Big( \br(t) \in \mathcal{B}_\epsilon(\bx) \Big) = \int_{\overline{\mathcal{B}}_\epsilon(\bx)} f_\Rconf(t,\bx) \hspace{0.75mm} \textnormal{d}\bx
\end{equation}

%\todo[inline]{Define the closed ball}
where, $\overline{\mathcal{B}}_\epsilon(\bx) = \{\xi \in \mathbb{R}^d : \hspace{2mm} \| \xi - \bx\| \leq \epsilon \}$ is the closed neighborhood, and  $f_\Rconf(t,\bx)$ is the probability density function (PDF) that defines the distribution of the swarm over the space $\Rconf$ and it has the following properties,

\begin{subequations} \label{eqn: probProp}
	\begin{align}
	f_\Rconf (t,\bx) &\geq \hspace{0.5mm} 0 \qquad \forall t \in \mathbb{R}^+, \ \forall \bx \in \Rconf, \\
	\hspace{-0.6mm}\int_{\Rconf} f_\Rconf (t,\bx) \hspace{0.5mm} \textnormal{d}\bx \hspace{0.6mm} &= \hspace{0.6mm} 1 \qquad \forall t \in \mathbb{R}^+.
	\end{align}
\end{subequations}

Here we consider a swarm with $N$ agents at time $t$ as a realization of the PDF $f_\Rconf (t,\bx)$, i.e., $N$ samples drawn from this PDF.   In the context of swarm density control, we will also refer to  $f_\Rconf (t,\bx)$ as the \emph{swarm density distribution}.

\begin{rem}
	The swarm density distribution definition is given with the assumption that the swarm is homogeneous i.e. all agents are identical. Although it is beyond the scope of this paper, presented framework can be generalized to swarms with heterogeneous agents using the notion of \emph{weighted distributions} \cite{patil2002weighted}. This approach adjusts $\rho(t,\bx)$ in (\ref{eqn: trueDens}) by putting an emphasis on each agent via the assigned individual weights (e.g., mass of each agent). 
\end{rem}

\subsection{Modeling the Motion of Agents} \label{sec: modeling}
The notion of swarm density distribution is independent of $N$. For single agent case ($N = 1$), $f_\Rconf(t,\bx)$ is the probability of finding the agent in a $\epsilon$-neighborhood of any point $\bx \in \Rconf$ at time $t$ as $\epsilon \rightarrow 0$.  Motivated by this interpretation, the motion of an agent  can be modeled with the following stochastic differential equation (SDE),

\begin{equation} \label{eqn: SDE}
	\textnormal{d}\br (t) = \bv\big(t,\br(t)\big)\textnormal{d}t + \sigma\big(t,\br(t)\big)\textnormal{d}\textbf{W}(t),
\end{equation}

where $\bv(t,\bx)$ is the  velocity field that acts on the agent and $\textbf{W}(t)$ is continuous-time stochastic process called \emph{Wiener process} with the diffusion function $\sigma\big(t,\br(t)\big)$. Given that agents' motion over $\Rconf$ is governed by (\ref{eqn: SDE}), the probability density function $f_\Rconf(t,\bx)$ satisfies the following \emph{Fokker-Planck equation} \cite{risken1984fokker}:

\begin{equation}
	{\frac {\partial }{\partial t}}f_\Rconf(t,\bx)=-\nabla\big[\bv(t,\bx)f_\Rconf(t,\bx)\big]+\Delta\big[\Sigma(t,\bx)f_\Rconf(t,\bx)\big],
\end{equation}

where $\Sigma(t,\bx) \!=\! \sigma^2 (t,\bx) / 2$ is the diffusion coefficient. 
With our interpretation of the swarm being a realization of the PDF $f_\Rconf(t,\bx)$,   we can state that the swarm density distribution is  propagated by using the Fokker-Planck equation. Within the scope of this paper we will ignore stochastic forces acting on the agents, i.e., $\sigma (t,\bx) = 0$. Consequently, equations are simplified as follows:

\begin{subequations}
	\begin{align}
	\dot{\br} (t) &= \bv\big(t,\br(t)\big), \label{eqn: simplifiedEqns1} \\
	\dot{f}_\Rconf(t,\bx) &= - \nabla\big[\bv(t,\bx)f_\Rconf(t,\bx)\big]. \label{eqn: simplifiedEqns2}
	\end{align}
\end{subequations}
\vspace{-1mm}
\begin{rem}
	Observe that equation (\ref{eqn: simplifiedEqns1}) is the constraint-free kinematics of an object and the equation (\ref{eqn: simplifiedEqns2}) is the \emph{continuity equation} in physics that governs the transport of a conserved  quantity.
\end{rem}

\subsection{Problem Definition}
Let $\Swarm(t) = \{\br_1(t),\br_2(t),\dots,\br_N(t)\}$ be a homogeneous swarm of $N$  agents that are distributed over configuration space i.e. $\br_i(t) \in \Rconf \ \forall i$ with initial swarm density distribution $f_\Rconf(t_0,\bx)$. For a given continuous desired density distribution $f_\Rconf^d(\bx)$, the \emph{density control} problem  of  this paper is defined as synthesizing a local density feedback law based on \emph{velocity field} $\bv(t,\bx)$ such that, when followed by agents, the swarm density distribution satisfies the following condition,

\begin{equation}
	\lim\limits_{t \rightarrow \infty} f_\Rconf(t,\bx) = f_\Rconf^d(\bx).
\end{equation}

\begin{figure}[hbt!]
\centering
\includegraphics[width=3.1in]{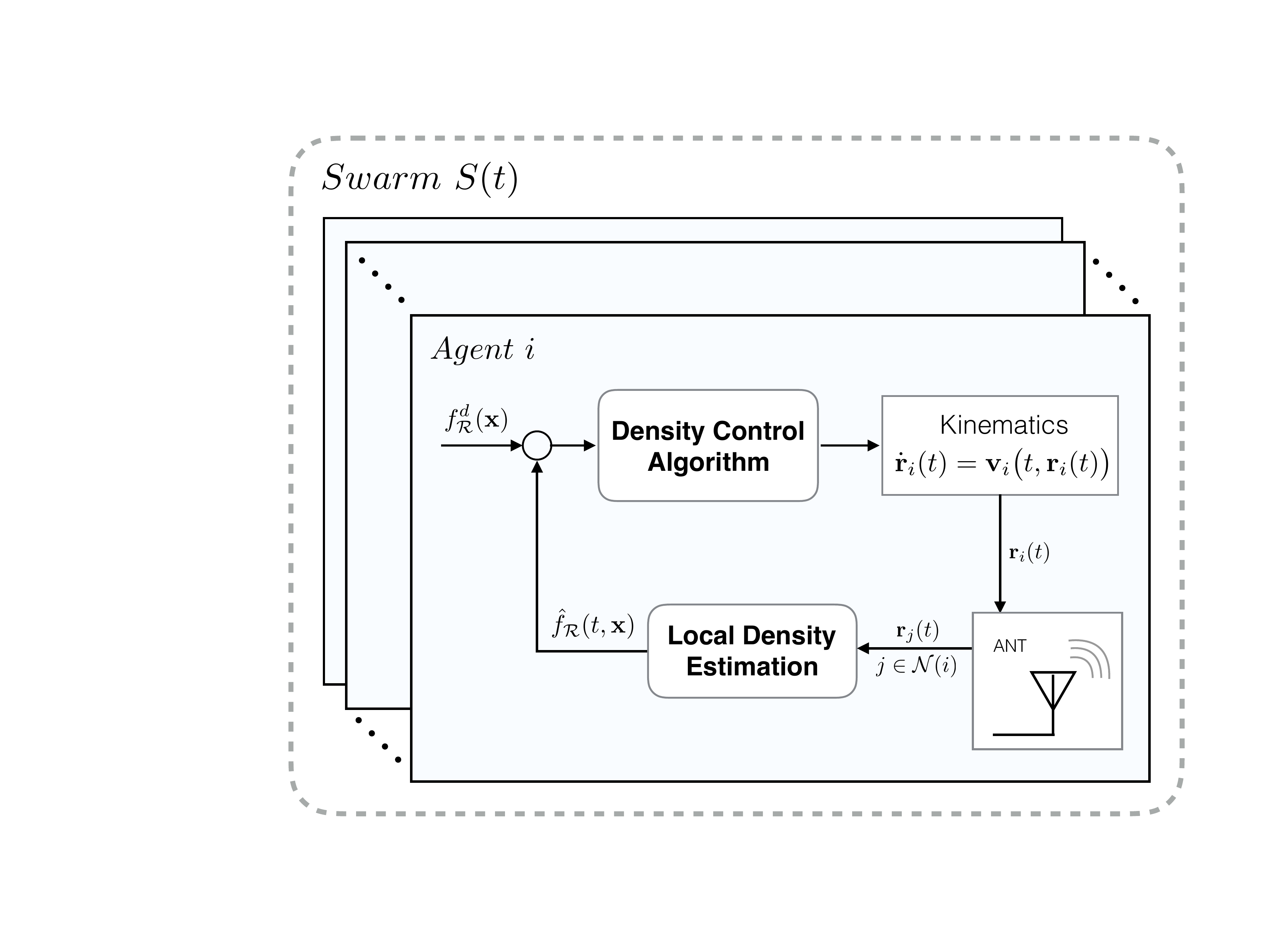}
\caption{The block diagram of the local density feedback control architecture for swarm of agents.} 
\label{fig: blockdiagfeedback}
\end{figure}

We  refer to the previous work \cite{swarm_denc13_ifac,swarm_denc13_acc,swarm_coll13,swarm_denc14_tac,demir2015probabilistic,demir2015decentralized} on how to design sequence of desired \emph{swarm density distributions} $f_\Rconf^d(\bx)$ with prescribed  properties (e.g., safety, motion and flow constraints), as it is not the  focus of this paper. As it is illustrated in Figure \ref{fig: blockdiagfeedback}, the idea is that, the desired density (analogous to reference signal in feedback control systems)  can be produced using the techniques presented in these earlier papers, and the density feedback law to be introduced   can be used to track them by generating a velocity field $\bv(t,\br)$ as a function of the desired density and the current density. 

%===============================================================================
%===============================================================================
\section{Local Density Estimation}%\label{}
The feedback  law to compute the velocity field requires  estimates of the \emph{local density distribution}, i.e., density distribution around each agent and its gradient.  To that end, before getting into the details of velocity field feedback law, this section addresses the local density estimation problem for the swarm. 

Having finite number of agents, $N$, prohibits the exact calculation of the local densities $\rho(t,\br_i(t))$, $\forall\br_i(t) \in S(t)$. Therefore, we aim to estimate the PDF, swarm density distribution, from given agent positions. For this purpose, we utilize   \emph{kernel density estimation} (KDE) \cite{rosenblatt1956remarks,parzen1962estimation}, which is devised to estimate the PDF of a random variable from a given sampling. The kernel density estimate of $f_\Rconf(t,\bx)$ for any point $\bx \in \Rconf$ at time $t$ is given by, 

\begin{equation}\label{eqn:est1}
	\hat{f}_\Rconf(t,\bx) = \bigintssss_{\Rconf} \ \Bigg[ \prod_{k=1}^{d} \frac{1}{h_k}K\bigg(\frac{x^{[k]}-\xi^{[k]}}{h_k}\bigg) \Bigg] \ \textnormal{d}P_N(t,\textbf{$\xi$})
\end{equation}

where $\textnormal{d}P_N(t,\textbf{$\xi$})$  is defined as,

\begin{equation*}
	\textnormal{d}P_N(t,\textbf{$\xi$}) = \frac{1}{N}\sum_{\br(t) \in S(t)} \delta(\textbf{$\xi$}-\br (t)) \ \textnormal{d}\xi,
\end{equation*}

where $\delta(\cdot)$ is the multi-dimensional Dirac $\delta$-function and substitution reduces (\ref{eqn:est1}) to the following form

\begin{equation} \label{eqn:f_estimator}
	\hat{f}_\Rconf (t,\bx) = \frac{1}{N} \sum_{i = 1}^{N} \Bigg[\prod_{k=1}^{d} \frac{1}{h_k}K\bigg(\frac{x^{[k]}-r^{[k]}_i(t)}{h_k}\bigg)\Bigg].
\end{equation}

In (\ref{eqn:f_estimator}),  $K(\cdot): \mathbb{R} \rightarrow \mathbb{R}$ is called \emph{smoothing kernel}  with \emph{smoothing parameter} $h_k$ in each dimension $k \in \{1,\dots, d\}$. Generally, smoothing parameters are chosen to be identical in each dimension, i.e., $h_k = h,\ \forall k$. Consequently, the local density distribution around an agent $i$ is estimated as follows,

\begin{equation} \label{eqn:estimator}
	\hat{f}_\Rconf (t,\br_i(t)) = \frac{1}{N h^d} \sum_{j = 1}^{N} \Bigg[\prod_{k=1}^{d} K\bigg(\frac{r^{[k]}_i(t)-r^{[k]}_j(t)}{h}\bigg)\Bigg].
\end{equation}

\subsection{Smoothing Kernel Function}
Kernel density estimation boils down to a proper choice  of kernel function $K(\cdot)$ and smoothing parameter $h$. Typically, kernel functions are chosen such that they satisfy  the following properties \cite{parzen1962estimation},

\begin{subequations}\label{eqn:Kconditions}
	\begin{align} 
		\sup_{-\infty < x < \infty} \big|K(x)\big| &< \infty, \\
		\int_{-\infty}^{\infty} \big|K(x)\big| \ dx &< \infty, \\
		\lim\limits_{x \rightarrow \infty} \big|xK(x)\big| &= 0.
	\end{align}	
\end{subequations}

The following theorem \cite{bochner1955harmonic} states the conditions for the PDF estimation (\ref{eqn:estimator}) to be asymptotically unbiased in the sense that if smoothing parameter $h = h(N)$ is chosen as function of $N$ such that  

\begin{equation} \label{eqn:hformula}
	\lim\limits_{N \rightarrow \infty} h(N) = 0.
\end{equation}
\vspace{1mm}

\begin{thm}
	Suppose $K(x)$ is a \emph{Borel function} satisfying the conditions (\ref{eqn:Kconditions}). Let $g(x)$ satisfy $\int_{-\infty}^{\infty} |g(x)| \ dx < \infty$ and let $h(N)$ be a sequence of positive constants satisfying (\ref{eqn:hformula}). Define
	
	\begin{equation} \label{eqn:proof1def}
		\hat{g}(x) = \frac{1}{h(N)} \int_{-\infty}^{\infty} K\bigg(\frac{\xi}{h(N)}\bigg) \ g(x - \xi) \ \textnormal{d}\xi,
	\end{equation}
	
	then the following relation holds for every point $x$
	
	\begin{equation} \label{eqn:proof1res}
		\lim\limits_{N \rightarrow \infty} \hat{g}(x) = g(x) \ \int_{-\infty}^{\infty} K(\xi) \ \textnormal{d} \xi.
	\end{equation}
	
	Consequently, the estimates defined by (\ref{eqn:estimator}) are \emph{asymptotically unbiased} at all points $\bx$ at which $f_\Rconf(t,\bx)$ is continuous if the function $K(x)$ satisfies (\ref{eqn:Kconditions}) and if $h$ satisfies (\ref{eqn:hformula}) along with the following,
	
	\begin{equation} \label{eqn:normalKernel}
		\int_{-\infty}^{\infty} K(\xi) \ \textnormal{d} \xi = 1.
	\end{equation}
\end{thm}

Typically  $K(x)$ is chosen to be  radially symmetric, and it is  unimodal to represent the PDF around a point, i.e., $K(x)$ has a single maximum that occurs at $x = 0$. Examples of such smoothing kernel functions are listed in \cite{parzen1962estimation}. Also note that the definition (\ref{eqn:proof1def})  with (\ref{eqn:hformula}),(\ref{eqn:proof1res}) and (\ref{eqn:normalKernel}), implies that 

\begin{equation*}
	\lim\limits_{N \rightarrow \infty} \Bigg(\frac{1}{h(N)} K\bigg(\frac{\xi}{h(N)}\bigg)\Bigg) = \delta(\xi).
\end{equation*}

In other words, as $h(N) \hspace{-0.7mm} \rightarrow \hspace{-0.7mm} 0$ the smoothing kernel approaches to Dirac $\delta$-function, hence, the estimation $ \hat{f}_\Rconf (t,\bx) $ approaches to the continuous actual swarm density distribution $f_\Rconf(t,\bx)$ as given by the Theorem 1 in \cite{wied2012consistency}:  

\begin{thm}[Weak consistency] \label{thm: probabilisticConvergence}
	Let the assumption
	
	\begin{equation*}
		\lim_{N \rightarrow \infty} N h(N) = \infty
	\end{equation*}
	
	hold. Then, at each point of continuity $x$ of $f$, the estimator $f_N(x)$ is weakly consistent, i.e. for each $\epsilon > 0$
	
	\begin{equation}
		\lim_{N \rightarrow \infty} \mathbb{P} (|f_N(x) - f(x)| > \epsilon) = 0
	\end{equation}
\end{thm}

%\begin{equation} \label{eqn: probabilisticConvergence}
%	\lim\limits_{N \rightarrow \infty} \textnormal{prob}\Big( \big|\hat{f}_\Rconf(t,\bx) - f_\Rconf(t,\bx)\big| > \epsilon \Big) = 0,
%\end{equation}
%
%when the condition $\lim_{N \rightarrow \infty} N h(N) = \infty$ is satisfied. 
Note that, for finite $N$, smoothing parameter $h$ can not be arbitrary small, as after certain $h_{min}$ the estimate $\hat{f}_\Rconf(t,\bx)$ becomes nothing but $N$ number of impulses at agent locations. 

\subsection{Selection of Smoothing Parameter}
The necessity of $ h > h_{min}$ (or finite $N$) incurs bias error for the estimation $\hat{f}_\Rconf (t,\bx)$. Also, selection of $h$ affects the behavior of the estimator in the sense that, if $h$ is small then the estimator will give results with high noise and variance, similarly if $h$ is large then the results will be too smooth, unable to capture the characteristics of the actual density. Typically, smoothing parameter $h$ is   chosen such that it minimizes the overall integral of mean-squared error (MSE) between $\hat{f}_\Rconf (t,\bx)$ and $f_\Rconf(t,\bx)$ which is defined as,

\begin{align} \label{eqn: biasVarTrade}
	\begin{split}
	\hspace{-0.17cm}\mathbb{E}\Big[\big(\hat{f}_\Rconf(t,\bx) - f_\Rconf(t,\bx) \big)^2 \Big] \hspace{-0.1cm} \ = \ & \textnormal{bias}\big(\hat{f}_\Rconf(t,\bx)\big)^2 \hspace{-0.05cm} \\ &+ \textnormal{var}\big(\hat{f}_\Rconf(t,\bx)\big)
	\end{split}
\end{align}

Note that the expression (\ref{eqn: biasVarTrade}) demonstrates the \emph{bias-variance trade-off} where the bias and variance expressions are as follows \cite{hansen2009lecture},

\begin{align} \label{eqn: bias}
	\textnormal{bias}\big(\hat{f}_\Rconf(t,\bx)\big) % &=  \mathbb{E}\Big[\hat{f}_\Rconf(t,\bx) - f_\Rconf(t,\bx) \Big] \\
										     &=  h^\nu \frac{\kappa_\nu}{\nu!} \sum_{k=1}^d \frac{\partial^\nu}{\partial x_k^\nu} f_\Rconf(t,\bx) + o(h^\nu),
\end{align}

\begin{align}  \label{eqn: variance}
	\hspace{-5mm}\textnormal{var}\big(\hat{f}_\Rconf(t,\bx)\big) %&= \mathbb{E}\big[\hat{f}_\Rconf(t,\bx)^2\big] - \mathbb{E}\big[\hat{f_\Rconf}(t,\bx)\big]^2  \\
										     &= \frac{f_\Rconf(t,\bx) R(K)^d}{N h^d} + O\Big(\frac{1}{N}\Big), 
\end{align}

where $\nu$ is the \emph{order} of $K(x)$ which is the index of first non-zero \emph{moment} $\kappa_j$:

\begin{equation*}
	\kappa_j = \int_{-\infty}^{\infty} x^j K(x) \ \textnormal{d}x ,
\end{equation*}

and $R(K)$ is called \emph{roughness} of $K(x)$ and calculated as,

\begin{equation*}
	R(K) = \int_{-\infty}^{\infty} K(x)^2 \ \textnormal{d}x.
\end{equation*}

$R(K)$ is used  as the measure of difficulty to estimate $K(x)$, i.e., larger the value of $R(K)$ more difficult it is to estimate $K(x)$ \cite{sheather2004density}.  Unfortunately, closed-form solution for $h$ that achieves the minimum for the integral of (\ref{eqn: biasVarTrade}) does not exists \cite{hansen2009lecture}. A rule of thumb for determining $h$ is to replace $f_\Rconf(t,\bx)$ in (\ref{eqn: bias}) and (\ref{eqn: variance}) with another elementary function (e.g., Gaussian) and find the minimizing $h$ as if the actual density $f_\Rconf(t,\bx)$ is the elementary function. The general form for minimizing $h$ is given as follows, 

\begin{equation}
	h^* = \hat{\sigma} \ C_\nu(K,d) \ N^{-1/(2\nu + d)},
\end{equation} 

where $\hat{\sigma}$ is the standard deviation of the $N$ sample, and $C_\nu(K,d)$ is a constant depends on the order of $K$, number of dimensions $d$ and the roughness of the elementary function that is assumed to be the actual density. Numerical $C_\nu(K,d)$ values for some elementary functions are given in Table 7 in \cite{hansen2009lecture}.

\begin{rem}
	The kernel function $K(x)$ determines the characteristics of the spikes in the estimation $\hat{f}_\Rconf(t,\bx)$ and this estimate is differentiable  provided that the kernel function is differentiable.
\end{rem}

\begin{rem} \label{rmk: error}
	As it can be seen from (\ref{eqn: bias}), for symmetric non-negative kernels i.e. \emph{second-order} kernels ($\nu = 2$), the bias error is $O(h^2)$. However, as it is mentioned in \cite{monaghan2005smoothed}, this error will be smaller when particles being propagated by certain dynamical equations that drive the particles  to  lower potential energy states. 
\end{rem}

%===============================================================================
%===============================================================================
\section{Velocity Field Generation}%\label{}
Given initial $f_\Rconf(t_0,\bx)$ and desired swarm density distributions $f_\Rconf^d(\bx)$, our previous approach \cite{swarm_denc13_ifac,swarm_coll13,demir2015probabilistic,demir2015decentralized} generates individual velocities for each agent in a probabilistic manner. %, and resulting motion from the probabilistic approach mimics the behavior of gas particles while satisfying \emph{Liouville's} theorem \cite{liouville1838note}. 
However, probabilistic approach has two significant challenges: $(i)$ collisions due to randomized motion, $(ii)$ ongoing motion of the agents even when $f_\Rconf(t,\bx) = f_\Rconf^d(\bx)$.
To overcome these challenges, we propose a velocity field computation method that uses the difference between the local density estimate and the desired density (i.e., the density error estimate). The resulting velocity field will be a continuous function in space, and hence, will significantly reduce collisions between agents. Furthermore, velocity field will dissipate as current density converges to the desired one.  In the derivation of velocity field synthesis method, the heat equation plays a critical role, which is introduced next.

\subsection{Velocity Field Generation using Heat Equation}
Heat equation is a partial differential equation that describes the evolution of temperature in time over a region:

\begin{equation}
	\frac{\partial u(t,\bx)}{\partial t} = D \hspace{0.5mm} \Delta u(t,\bx),
\end{equation}

where $D > 0$ and $\Delta$ is Laplacian operator describing the sum of second spatial derivatives as follows,

\begin{equation*}
	\Delta(\cdot) = \sum_{i=1}^d \frac{\partial^2 (\cdot)}{\partial x^2_i} \ .
\end{equation*}

This equation is also known as a special form of \emph{diffusion equation} in which the \emph{diffusion constant} $D$ is considered to be a function of $\bx$ and $u$. In this paper, heat equation will be used to describe the evolution of the swarm density distribution, PDF given by $f_\Rconf(t,\bx)$, in time. 

Let the difference between the  current and desired probability density functions be 
\vspace{2mm}
\begin{equation*}
	\Phi(t,\bx) := f_\Rconf(t,\bx) - f_\Rconf^d(\bx)
\end{equation*} 

over $\Rconf$.
The core idea of this paper is to find the velocity field that transforms Equation (\ref{eqn: simplifiedEqns2}) to heat equation given in (\ref{eqn: heatEq1}). For this purpose, we propose the following feedback law to compute the velocity as a function of $\Phi(t,\bx)$ and $f_\Rconf(t,\bx)$, 

\begin{equation}
	\label{eqn: controlLaw}
	\bv(t,\bx) = -D \ \frac{\nabla\Phi(t,\bx)}{f_\Rconf(t,\bx)} .
\end{equation} 

Noting that the desired density is constant in time, i.e. $\partial f_\Rconf^d(\bx) / \partial t = 0$, we can rewrite (\ref{eqn: simplifiedEqns2}) with the proposed velocity field as follows,

\begin{subequations} \label{eqn: heatEq}
	\begin{align}
		\frac{\partial \Phi} {\partial t} (t,\bx) &= D \hspace{0.5mm} \Delta \Phi (t,\bx), \quad \bx \in \Rconf, \ t > t_0 \label{eqn: heatEq1} \\
		\Phi(t_0,\bx) &= f_\Rconf(t_0,\bx) - f_\Rconf^d(\bx).
	\end{align}
\end{subequations}

and suppose that $\Phi(t,\bx)$ satisfies the boundary condition $\nabla\Phi(t,\bx) = 0$ on $\partial\Rconf$. Then $\int_\Rconf \Phi(t,\bx) \hspace{0.5mm} \textnormal{d}\bx $ is conserved for all $t \geq t_0$, and, % Recall that density estimates are normalized with the total number of agents as it is stated in (\ref{eqn:estimator}). 
since we have,  %a swarm where the total number of agents is conserved over $\Rconf$ i.e.,

\begin{equation} \label{eqn: densEqual}
	\int_\Rconf f_\Rconf(t,\bx) \hspace{0.5mm} \textnormal{d}\bx = \int_\Rconf f_\Rconf^d(\bx) \hspace{0.5mm} \textnormal{d}\bx \quad \forall t \geq t_0
\end{equation}

the following holds,

\begin{equation} \label{eqn: IntPhi}
	\int_\Rconf \Phi(t,\bx) \ \textnormal{d}\bx = 0 \quad \forall t\geq t_0.
\end{equation}

%Also, it is crucial to note that the solution $\Phi(t,\bx)$ for the problem (\ref{eqn: heatEq}) with \emph{Neumann} boundary conditions is \emph{unique}  \cite{barton1989elements}. 

The following theorem gives the main result which shows that, in case of $N \rightarrow \infty$, the swarm converges to the desired density distribution from any initial distribution when the velocity field is synthesized with (\ref{eqn: controlLaw}). Also recall from Theorem \ref{thm: probabilisticConvergence} that as $N \rightarrow \infty$, with probability $1$, $\hat{f}_\Rconf \rightarrow f_\Rconf$. \\

\begin{thm}
	Consider a swarm $\Swarm(t)$, commanded with continuous desired density function $f_{\Rconf}^d(\bx)$, where the motion of each agent is governed by the following  dynamics,
	
	\begin{equation} \label{eqn: agentDyn}
		\dot{\bx}(t) = -D \ \frac{\nabla \Phi(t,\bx(t))}{f_\Rconf(t,\bx(t))}.
	\end{equation}
	
	As $t \rightarrow \infty$, for any initial conditions of the agents in the swarm, the density difference $\Phi(t,\bx(t))$ for each agent asymptotically converges to zero i.e.
	
	\begin{equation*}
		\lim\limits_{t \rightarrow \infty} \Phi(t,\bx) = 0 \quad \bx \in \Rconf,
	\end{equation*}
	
	hence the swarm density distribution  $f_\Rconf(t,\bx(t))$ approaches to the desired density distribution  $f_\Rconf^d(\bx(t))$ and velocities of all the agents vanish asymptotically,  i.e., $\lim_{t \rightarrow \infty} \dot{\bx}(t) = 0$.
\end{thm}

\begin{pf}
	For the analysis we define the following positive definite function as the Lyapunov function for the swarm:
	
	\begin{align}
		V(t) &= \frac{1}{2} \int_\Rconf \bigg(\frac{f_\Rconf(t,\bx)}{D}\bigg)^2 \dot{\bx}^T \dot{\bx} \ \textnormal{d} \bx \\
		&= \frac{1}{2} \int_\Rconf \nabla\Phi(t,\bx)^T \nabla\Phi(t,\bx) \ \textnormal{d} \bx.
	\end{align}
	
	 Taking the time derivative of Lyapunov function $V(t)$ yields,
	
	\begin{equation*}
		\dot{V}(t) = \int_\Rconf \nabla\Phi(t,\bx)^T \nabla\dot{\Phi}(t,\bx) \ \textnormal{d} \bx,
	\end{equation*}
	
	by using heat equation $\dot{\Phi}(t,\bx) = D \hspace{0.3mm} \Delta\Phi(t,\bx)$ and the third derivative vector calculus,
	
	\begin{align*}
		\dot{V}(t)	&= \int_\Rconf \nabla\Phi(t,\bx)^T \Big(\nabla\big(D \hspace{0.3mm} \Delta\Phi(t,\bx)\big) \Big) \textnormal{d} \bx \\
		&= D \int_\Rconf \nabla\Phi(t,\bx)^T \Delta\big(\nabla\Phi(t,\bx)\big) \ \textnormal{d} \bx.
	\end{align*}
Employing a change of variables as $\zeta(t,\bx) = \nabla\Phi(t,\bx)$ yields
	\begin{align*}
		\dot{V}(t) = D \int_\Rconf  \zeta(t,\bx)^T\Delta\zeta(t,\bx) \ \textnormal{d} \bx \\
	 \end{align*}
	 Also the boundary condition $\nabla\Phi(t,\bx) = 0$ on $\partial \Rconf$ can be written as \emph{Dirichlet} boundary condition of  $\zeta(t,\bx) = 0$ on $\partial \Rconf$. Then we can conclude $\dot{V}(t) < 0$ because of the fact that, \emph{Dirichlet} eigenvalue problem for $\Delta$ operator has countably many strictly negative eigenvalues  \cite{li1983schrodinger}. Consequently, $\lim_{t \rightarrow \infty} \nabla\Phi(t,\bx(t)) = 0$, which implies $\lim_{t \rightarrow \infty} \dot{\bx}(t) = 0$. Furthermore, 
\vspace{2mm}
	\begin{equation} \label{eqn: biasedDensity}
		\lim\limits_{t \rightarrow \infty} \Phi(t,\bx(t)) = \textnormal{const}.
	\end{equation}
	Since $f_\Rconf^d(\bx)$ is continuous then $\Phi(t,\bx)$ is also continuous and because of (\ref{eqn: IntPhi}), %
    $\lim\limits_{t \rightarrow \infty} \Phi(t,\bx(t)) = 0$ which implies 
\vspace{2mm}
	\begin{equation*}
		\lim\limits_{t \rightarrow \infty} f_\Rconf(t,\bx(t)) = f_\Rconf^d(\bx(t)). 
	\end{equation*}
	\hspace{8.25cm} \qed
\end{pf}
\begin{rem}
	The result (\ref{eqn: biasedDensity}) implies that even if the relation (\ref{eqn: densEqual}) does not hold, the desired density $f_\Rconf^d(\bx)$ will be achieved with a bias of $\Phi(t,\bx) = \textnormal{const}$. This means if the amount of agents over the domain $\Rconf$ increase or decrease at an instant without tuning $N$ in (\ref{eqn:estimator}), the velocity field will act to preserve the characteristics of the desired density $f_\Rconf^d(\bx)$.
\end{rem}

The velocity field requires the local density distribution and its gradient, and this information is estimated on-board of each agent with the local information, hence, making the algorithm decentralized. By each agent following this local velocity accurately with the feedback of local density estimate, the swarm will be driven to the desired density $f_\Rconf^d(\bx)$ autonomously. Also, as it is mentioned in Remark \ref{rmk: error}, since the heat equation, and hence, the velocity field spreads agents locally uniformly, the bias error will be better than $O(h^2)$. Next section provides simulation results justifying these claims.

%===============================================================================
%===============================================================================
\section{Simulation Results}\label{sims}

\begin{figure}[hbt!]
\hspace{-5mm}\vspace{-1.5cm}\includegraphics[width=3.7in]{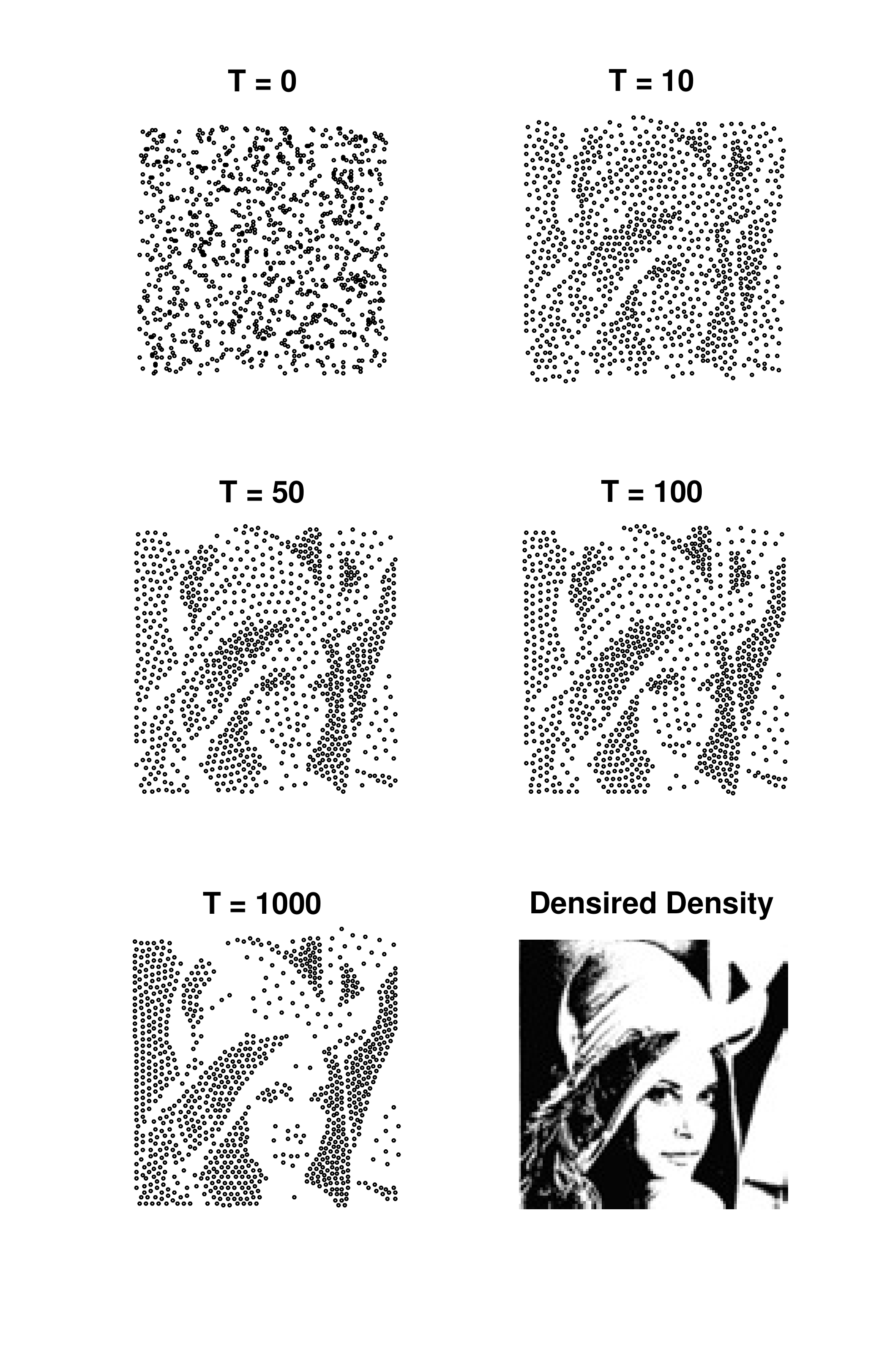}
\caption{The swarm density distribution in time towards desired density.} 
\label{fig: pattern}
\end{figure}

In this simulation, we illustrate the capability of the method to generate complex patterns. In Figure \ref{fig: pattern}, a swarm of $N = 1000$ agents starting from a random sample of positions from uniform distribution are commanded to desired density defined by the picture of Lenna. The picture file processed to construct a normalized continuous desired swarm density function. The local density estimator utilizes a multi-dimensional Gaussian kernel function:
\vspace{3mm}
\begin{equation*}
	K(\bx) = \frac{2}{\pi} \exp\bigg({-2\frac{\bx^T\bx}{h^2}}\bigg), 
\end{equation*}

and the smoothing parameter $h$ is taken as $h = L/20$ where $L$ is a dimension of a square domain. The effective radius around each agent is $R_{eff} = 2h$. The heat equation based controller has the diffusion constant of $D = 5$. Results show that swarm of 1000 agents almost converged to the desired density distribution after $T = 1000$. Also, agents are locally uniformly distributed and inter-agent distances are determined by the density command. 

\begin{figure}[hbt!]
\centering
\includegraphics[width=3.2in]{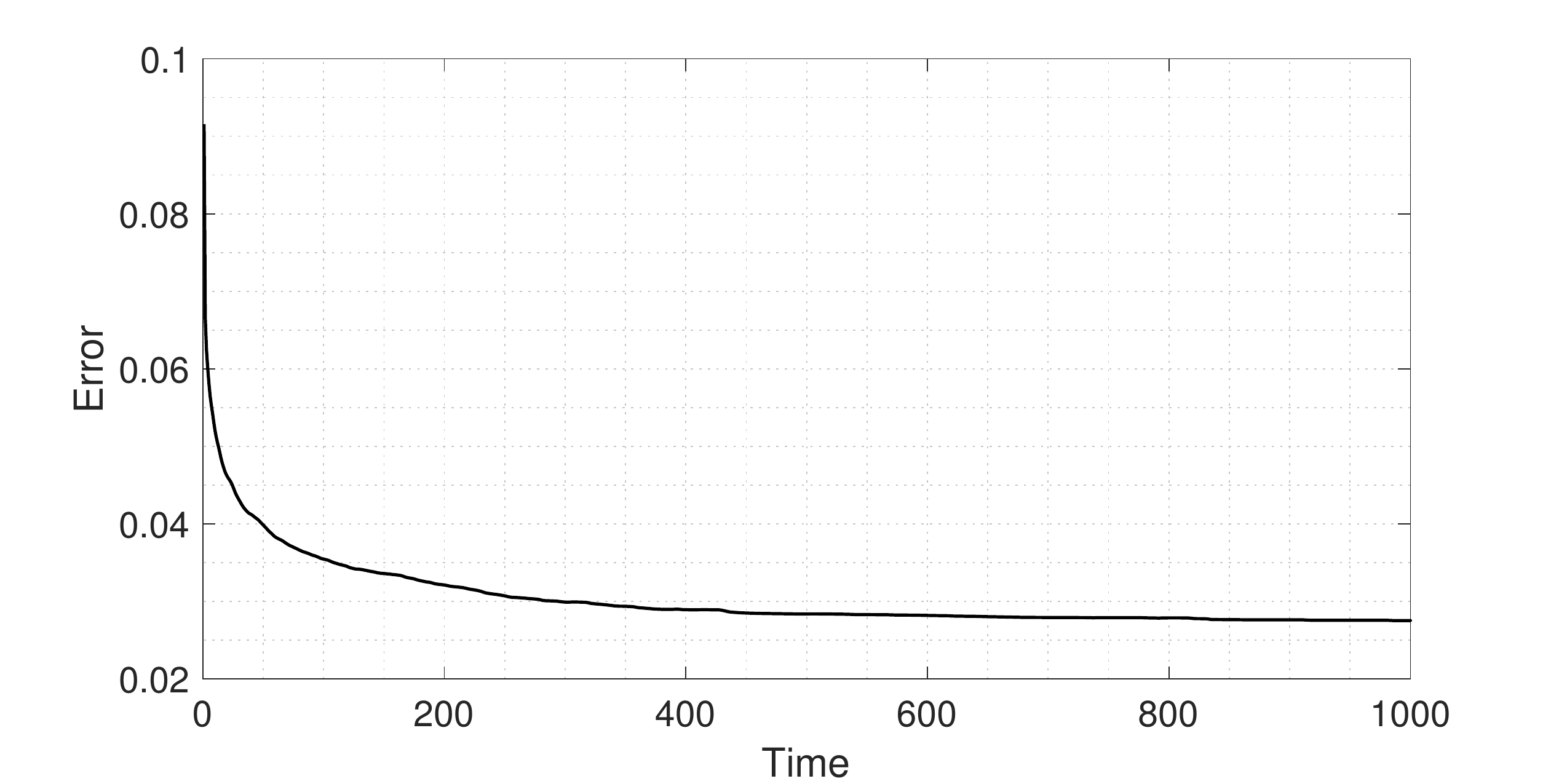}
\caption{The MSE error in time from desired distribution.} 
\label{fig: MSE}
\end{figure}

Figure \ref{fig: MSE} shows the integral of point-wise error over the domain in time defined as,

\begin{equation*}
	E(t) = \int_\Rconf \big|\hat{f}_\Rconf(t,\bx) - f_\Rconf^d(\bx)\big| \hspace{0.5mm} \textnormal{d}\bx
\end{equation*}

Even though the error dissipates in time, a permanent error exists due to the value of smoothing parameter $h$, i.e., estimation error from the kernel density estimation.

%===============================================================================
%===============================================================================
\section{Conclusion}
In this paper, we have presented a velocity field synthesis method for a swarm of agents controlled with density distributions over a bounded domain. The synthesis method is decentralized in the sense that it is based on local density distribution feedback, which is achieved with kernel density estimation and utilized in a heat equation based control law. The key advantage of this method is that the resulting velocity field is smooth and it facilitates collision avoidance. We have also provided asymptotic convergence and stability analysis of the proposed velocity field for the case $N \rightarrow \infty$. In future work, we will analyze the method for a finite number of agents with estimation errors, we will also explore the stochastic part of the Fokker-Planck equation to enrich modeling capabilities. Finally, we will provide analysis for collision avoidance and expand the method to more complex dynamics with force fields. 

%\begin{ack}
%This research was supported in part by Defense Advanced Research Projects Agency (DARPA) Grant No. %D14AP00084, and by National Science Foundation (NSF) Grant No. CNS-1624328.
%\end{ack}

\bibliography{root}             % bib file to produce the bibliography
                                                     % with bibtex (preferred)
                                                   
%\appendix
%\section{A summary of Latin grammar}    % Each appendix must have a short title.
%\section{Some Latin vocabulary}              % Sections and subsections are supported  
                                                                         % in the appendices.
\end{document}